\documentclass[portuges,12pt,letter]{article}
\usepackage[centertags]{amsmath}
\usepackage{amsfonts}
\usepackage{newlfont}
\usepackage{amscd}
\usepackage{graphics}
\usepackage{epsfig}
\usepackage{indentfirst}
\usepackage{amsxtra}
\usepackage[latin1]{inputenc}
\usepackage{amssymb, amsmath}
\usepackage{amsthm}
\usepackage[mathscr]{eucal}

\title{\bf  An approximate numerical method for ordinary differential equation systems with applications to a flight mechanics model}

\author{Fabio Silva Botelho \\  Department of Mathematics \\
Federal University of Santa Catarina \\
Florian\'{o}polis - SC, Brazil}
\date{}
\begin{document}
\maketitle

\abstract{This short communication develops a new numerical procedure suitable for a large class of ordinary differential equation systems found in models in physics and engineering.
The main numerical procedure is analogous to those concerning the generalized method of lines, originally published  in the here referenced books of 2011 and 2014, \cite{901,12} respectively. Finally, in the last section, we apply the method to a model in flight mechanics.  }

\section{Introduction}
Consider the first order system of ordinary differential equations given by
$$ \frac{ d u_j}{dt}=f_j(\{u_l\}), \text{ on } [0,t_f]\; \forall j \in \{1,\cdots,4\},$$ with the boundary conditions
$$u_1(0)=0,\; u_2(0),\; u_4(0)=0,\; u_4(t_f)=u_f.$$

Here $\mathbf{u}=\{u_l\} \in W^{1,2}([0,t_f],\mathbb{R}^4)$ and $f_j$ are functions at least of $C^1$ class, $\forall j \in \{1,2,3,4\}.$

Our proposed method is iterative so that we choose an starting solution denoted by $\tilde{u}$.

At this point, we define the number of nodes on $[0,t_f]$ by $N$ and set $d=t_f/N.$

Similarly to a proximal approach, we propose the following algorithm (in a similar fashion as those found in \cite{901,12,700}).

\begin{enumerate}
\item Choose $\tilde{u} \in W^{1,2}([0,t_f],\mathbb{R}^4).$
\item\label{AR1} Solve the equation system

\begin{equation}\label{AR2} K\frac{ d u_j}{dt}-(K-1)\frac{ d \tilde{u}_j}{dt}=f_j(\{u_l\}), \text{ on } [0,t_f]\; \forall j \in \{1,\cdots,4\},\end{equation} with the boundary conditions
$$u_1(0)=0,\; u_2(0),\; u_4(0)=0,\; u_4(t_f)=u_f.$$
\item Replace $\tilde{u}$ by $u$ and go to item (\ref{AR1}) up to the satisfaction of an appropriate convergence criterion.
\end{enumerate}

In finite differences, observe that the system may be approximated by
$$(u_j)_n-(u_j)_{n-1}=\frac{K-1}{K}((\tilde{u}_j)_n-(\tilde{u}_j)_{n-1})+f_j(\{(u_l)_{n-1}\})\frac{d}{K}$$

In particular, for $n=1$, we get
$$(u_j)_1=(u_j)_0+f_j(\{(u_l)_1\})\frac{ d }{K}+\frac{K-1}{K}((\tilde{u}_j)_1-(\tilde{u}_j)_0)+(E_j)_1,$$
where
$$(E_j)_1=[f_j(\{(u_l)_0\})-f_j(\{(u_l)_1\})]\frac{d}{K} \approx \mathcal{O}\left(\frac{d^2}{K}\right).$$

Similarly, for $n=2$ we have
\begin{eqnarray}
(u_j)_2-(u_j)_0&=& [(u_j)_2-(u_j)_1]+[(u_j)_1-(u_j)_0] \nonumber \\ &=&
f_j(\{(u_l)_1\})\frac{ d }{K}+\frac{K-1}{K}((\tilde{u}_j)_2-(\tilde{u}_j)_1)
\nonumber \\ &&+f_j(\{(u_l)_1\})\frac{ d }{K}+\frac{K-1}{K}((\tilde{u}_j)_1-(\tilde{u}_j)_0)+(E_j)_1
\nonumber \\ &=& 2f_j(\{(u_l)_1\})\frac{ d }{K}-2f_j(\{(u_l)_2\})\frac{ d }{K}+2f_j(\{(u_l)_2\})\frac{d}{K}
\nonumber \\ && +\frac{K-1}{K}((\tilde{u}_j)_2-(\tilde{u}_j)_0)+(E_j)_1, \end{eqnarray}

Summarizing
$$(u_j)_2= (u_j)_0 + 2f_j(\{(u_l)_2\})\frac{ d }{K}+\frac{K-1}{K}((\tilde{u}_j)_2-(\tilde{u}_j)_0)+(E_j)_2,$$
where $$(E_j)_2=(E_j)_1+\left(f_j(\{(u_l)_1\})-f_j(\{(u_l)_2\})\right)\frac{2d}{K}.$$

Reasoning inductively, for all $1\leq k\leq N,$ we obtain
$$(u_j)_k=(u_j)_0+f_j(\{(u_l)_k\})\frac{k d}{K}+\frac{K-1}{K}((\tilde{u}_j)_k-(\tilde{u}_j)_0)+(E_j)_k,$$
where,
\begin{eqnarray}(E_j)_k&=& (E_j)_{k-1}+k(f_j(\{u_l\}_{k-1})-f_j(\{u_l\}_k)\frac{ d }{K} \nonumber \\ &=&\sum_{m=1}^k(k-m+1)[f_j(\{(u_l)_{m-1}\})-f_j(\{(u_l)_m\})]\frac{d}{K} \approx \mathcal{O}\left(\frac{(k^2+k) d^2}{2K}\right),\end{eqnarray}
$\forall j \in \{1, \cdots,4\}.$

In particular, for $k=N$, for $K>0$ sufficiently big, we obtain

\begin{eqnarray}\label{AR3}(u_j)_N&=&(u_j)_0+f_j(\{(u_l)_N\})\frac{N d}{K}+\frac{K-1}{K}((\tilde{u}_j)_N-(\tilde{u}_j)_0)+(E_j)_N \nonumber \\ &\approx&
(u_j)_0+f_j(\{(u_l)_N\})\frac{N d}{K}+\frac{K-1}{K}((\tilde{u}_j)_N-(\tilde{u}_j)_0),
\forall j \in \{1,2,3,4\}.
\end{eqnarray}

In such a system we have $4$ equations suitable to find the unknown variables $(u_0)_3,(u_N)_1,(u_N)_2,(u_N)_3$, considering that
$(u_N)_4=(u_4)_f$ is known.

Through the system (\ref{AR3}) we obtain $u_N$ through the Newton's Method for a system of only the $4$ indicated variables.

Having $u_N$, we obtain $u_{N-1}$ through the equations

$$(u_j)_{N}-(u_j)_{N-1}\approx f_j(\{(u_l)_{N-1}\})\frac{d}{K}+\frac{K-1}{K}((\tilde{u}_j)_N)-(\tilde{u}_j)_{N-1}),$$
also through the Newton's Method.

Similarly, having $u_{N-1}$ we obtain $u_{N-2}$ through the system
$$(u_j)_{N-1}-(u_j)_{N-2}\approx f_j(\{(u_l)_{N-2}\})\frac{d}{K}+\frac{K-1}{K}((\tilde{u}_j)_{N-1})-(\tilde{u}_j)_{N-2}),$$
and so on up to finding $u_1.$

Having calculated $u$, we replace $\tilde{u}$ by $u$ and repeat the process up to an appropriate convergence criterion is satisfied.

The problem is then approximately solved.

\section{Applications to a flight mechanics model}

We present numerical results for the following system of equations, which models the in plane climbing  motion of
an airplane (please, see more details in \cite{127}).

\begin{equation}\left\{
\begin{array}{l}
 \dot{h}=V\sin \gamma,
 \\
 \dot{\gamma}=\frac{1}{m_f V}(F\sin[a+a_F]+L)-\frac{g}{V} \cos\gamma, \\
 \dot{V}=\frac{1}{m_f}(F\cos[a+a_F])-D)-g\sin\gamma \\
 \dot{x}=V\cos \gamma, \end{array} \right. \end{equation}

 with the boundary conditions,
 \begin{equation}\left\{
\begin{array}{l}
 h(0)=h_0,
 \\
 V(0)=V_0 \\
 x(0)=x_0 \\
 h(t_f)=h_f, \end{array} \right. \end{equation}

 where $t_f=532s$, $h$ is the airplane altitude, $V$ is its speed, $\gamma$ is the angle between its velocity and
 the horizontal axis, and finally $x$ denotes the horizontal coordinate position.

 For numerical purposes, we assume (Air bus 320)

 $m_f=120,000 Kg,$ $S_f=260 m^2$, $a= 0.0945 \text{ rad}$,
 $g=9.8 m/s^2$, $$\rho(h)=1.225(1-0.0065h/288.15)^{4.225}Kg/m^3,$$
 $$a_F=0.03225,$$
 $$(C_L)_a=5$$
 $$(C_D)_0=0.0175,$$
 $$K_1=0.0,$$
 $$K_2=0.06$$
 $$C_L=(C_D)_0+K1 C_L+K_2C_L^2,$$
 $$L=\frac{1}{2}\rho(h) V^2C_L S_f,$$
 $$D=\frac{1}{2}\rho(h) V^2C_DS_f,$$
 $$F=240000$$
 and
 where units refer to the International System.

 To simplify the analysis, we redefine the variables as below indicated:

 \begin{equation}\left\{
\begin{array}{l}
 h=u_1,
 \\
 \gamma=u_2 \\
 V=u_3 \\
 x=u_4. \end{array} \right. \end{equation}

Thus, denoting $\mathbf{u}=(u_1,u_2,u_3,u_4) \in U=W^{1,2}([0,t_f];\mathbb{R}^4),$ the system above indicated may be expressed by

 \begin{equation}\label{a.6v}\left\{
\begin{array}{l}
 \dot{u}_1=f_1(\mathbf{u})
 \\
 \dot{u}_2=f_2(\mathbf{u})\\
 \dot{u}_3=f_3(\mathbf{u}) \\
 \dot{u}_4=f_4(\mathbf{u}), \end{array} \right. \end{equation}

where,

\begin{equation}\left\{
\begin{array}{l}
 f_1(\mathbf{u})=u_3\sin (u_2),
 \\
 f_2(\mathbf{u})=\frac{1}{m_f u_3}(F\sin[a+a_F]+L(\mathbf{u}))-\frac{g}{u_3} \cos (u_2), \\
 f_3(\mathbf{u})=\frac{1}{m_f}(F\cos[a+a_F]-D(\mathbf{u}))-g\sin (u_2) \\
f_4(\mathbf{u})= u_3 \cos (u_2). \end{array} \right. \end{equation}

We solve this last system for the following boundary conditions:
 \begin{equation}\left\{
\begin{array}{l}
 h(0)=0\; m,
 \\
 V(0)=150  m/s, \\
 x(0)=0\;m, \\
 h(t_f)=11000\;m. \end{array} \right. \end{equation}

\begin{figure}
\centering \includegraphics [width=3in]{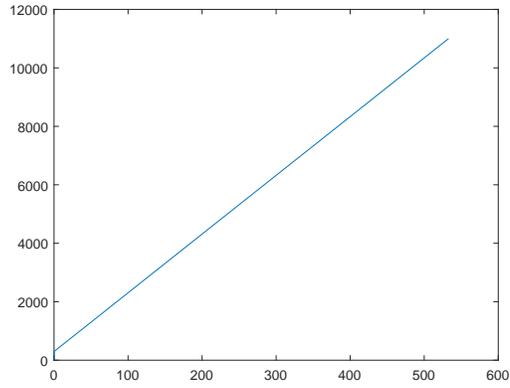}
\\ \caption{\small{The solution $h$ (in m) for $t_f=532s$.
}}\label{fig.1v}
\end{figure}

\begin{figure}
\centering \includegraphics [width=3in]{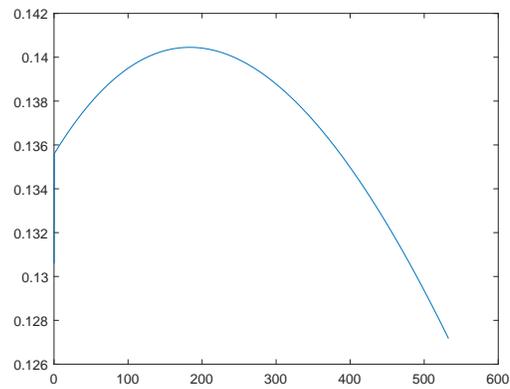}
\\ \caption{\small{The solution $\gamma$ (in rad) for $t_f=532s$.
}}\label{fig.2v}
\end{figure}

\begin{figure}
\centering \includegraphics [width=3in]{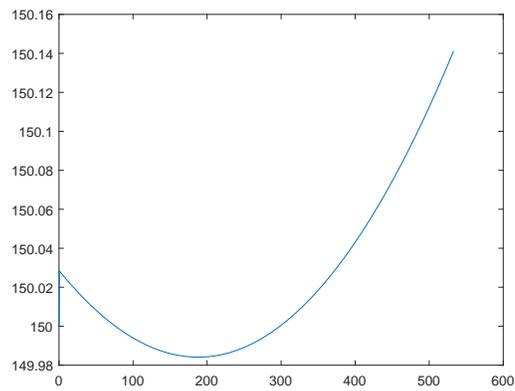}
\\ \caption{\small{The solution $V$ (in m/s) for $t_f=532s$.
}}\label{fig.3v}
\end{figure}

\begin{figure}
\centering \includegraphics [width=3in]{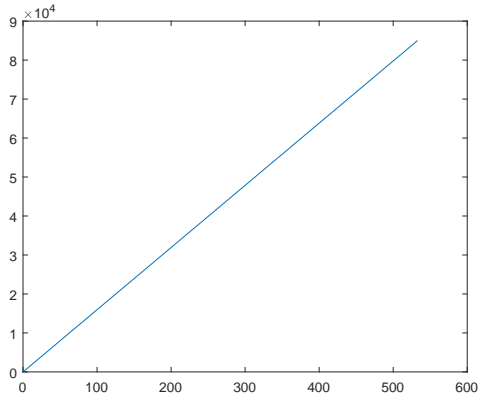}
\\ \caption{\small{The solution $x$ (in m) for $t_f=532s$.
}}\label{fig.4v}
\end{figure}

We have obtained the following solutions for $h,\gamma, V \text{ and } x$. Please see figures \ref{fig.1v}, \ref{fig.2v},  \ref{fig.3v} and \ref{fig.4v}, respectively. We have set $N=3000$ nodes and $K=509$.


\section{Acknowledgements}

The author is grateful to Professor Pedro Paglione of Technological Institute of Aeronautics, ITA, SP-Brazil, for his valuable suggestions and comments which help me a lot to improve some important parts
of this text, in particular on the part concerning the model in flight mechanics addressed.

\end{document}